\documentclass{article}
\usepackage{latexsym}
\usepackage{amssymb}
\begin{document}
\def\P{\Bbb P}
\def\Q{\Bbb Q}
\title{Higher-dimensional Forcing}
\author{BERNHARD IRRGANG}
\date{}
\maketitle
\begin{abstract}
We present a method of constructing ccc forcings: Suppose first that a continuous, commutative system of complete embeddings between countable forcings indexed along $\omega_1$ is given. Then its direct limit satisfies ccc by a well-known theorem on finite support iterations. However, this limit has size at most $\omega_1$. To get larger forcings, we do not consider linear systems but higher-dimensional systems which are indexed along simplified morasses.
\end{abstract}
\section{Introduction} 
In 1965, R. Solovay and S. Tennenbaum invented iterated forcing \cite{SolovayTennenbaum}. Since then it has become presumably the most important tool for proving relative consistency results in set theory. They developed it to prove the consistency of Suslin's hypothesis, which can be formulated as the non-existence of so-called Suslin trees. Given a Suslin tree, it is easy to destroy it by forcing. Unfortunately, in the generic extension, there can be new Suslin trees. They in turn can be destroyed by forcing and the right bookkeeping allows to finally destroy all Suslin trees. The main technical observation to make this work is that the generic extension of a generic extension is itself generic. More precisely, if $G$ is $\P_0$-generic over $M$ and $H$ is $\Q$-generic over $M[G]$, then there is a $\P_1$ such that $M[G][H]$ is $\P_1$-generic over $M$ and $\P_0$ completely embedds into $\P_1$. If we iterate this process, we obtain a commutative system $\langle \sigma_{\nu\tau}:\P_\nu \rightarrow \P_\tau \mid \nu \in \lambda\rangle$ of complete embeddings between forcings $\langle \P_\nu  \mid \nu \in \lambda\rangle$. To carry on this construction over limit ordinals, we can for exaple take the direct limit of the system at limit ordinals. This yields so-called finite support iterations. Taking inverse limits if $cf(\lambda)=\omega$ and dircet limits if $cf(\lambda)>\omega$ yields so-called countable support iterations. All this is nowadays well known to every set theorist.
\smallskip\\
Despite its huge success, there are important consistency questions that cannot be aswered by iterated forcing as we presently understand it. The most famous problem with iterated forcing is that some questions seem to require countable support iterations while it is nearly impossible to make $2^{\aleph_0}$ bigger than $\aleph_2$ with countable support iterations. An attempt to overcome such difficulties is given by M. Groszek and T. Jech \cite{JechGroszek}. They introduce iterations along arbitrary partially ordered sets instead only along ordinals.  
\smallskip\\
In the following, we persue a similar idea. We consider systems of embeddings between forcings that are indexed along simplified gap-1 or gap-2 morasses. However, not all embeddings will be complete embeddings. We use our approach to construct ccc forcings. Suppose first that a continuous, commutative system of complete embeddings between countable forcings indexed along $\omega_1$ is given. Then its direct limit satisfies ccc by a well-known theorem on finite support iterations. However, this limit has size at most $\omega_1$. To get larger forcings, we do not consider linear systems but higher-dimensional systems which are indexed along simplified morasses.
\smallskip\\
All our applications are of a rather combinatorial manner and we always use finite conditions. Therefore, only names for forcings play a role in our iterations while we do not need names for conditions. There are a number of results that can be reproved with our method.
P. Koszmider \cite{Koszmider} proved that it is consistent that there exists a sequence $\langle X_\alpha \mid \alpha < \omega_2\rangle$ of subsets $X_\alpha \subseteq \omega_1$ such that $X_\beta-X_\alpha$ is finite and $X_\alpha -X_\beta$ is uncountable for all $\beta < \alpha < \omega_2$. He uses S. Todorcevic's method of ordinal walks \cite{Todorcevic2007}. It is also known as the method of $\rho$-functions \cite{Koszmider1998} and provides a powerful tool to construct ccc forcings in the presence of $\Box_{\omega_1}$.  Other applications of it are a ccc forcing that adds an $\omega_2$-Suslin tree \cite{Todorcevic2007}, a ccc forcing for $\omega_2 \not\rightarrow (\omega:2)^2_\omega$ \cite{Todorcevic2007}, a ccc forcing to add a Kurepa tree \cite{Todorcevic2007, Velickovic} and a ccc forcing to add a thin-very tall superatomic Boolean algebra \cite{Todorcevic2007}. The last forcing was first found by Baumgartner and Shelah \cite{BaumgartnerShelah} independently from $\rho$-functions. That there can be a ccc forcing for $\omega_2 \not\rightarrow (\omega:2)^2_\omega$ was first observed by Galvin \cite{Kanamori}. That $\Box_{\omega_1}$ implies the existence of a ccc forcing which adds a Kurepa tree was first proved by Jensen \cite{Jensen2,Jensen3}. As we will see, all these results can be reproved from a simplified $(\omega_1,1)$-morass instead of $\Box_{\omega_1}$. The proof of Koszmider's result was already given in \cite{Irrgang5}.
\smallskip\\
All these examples add a structure on $\omega_2$. The natural question arises if something similar can be done for forcings that add a structure on a higher cardinal. This seems in general possible. As an example, we showed \cite{Irrgang6} that if there exists a simplified $(\omega_1,2)$-morass, then there exists a ccc forcing which adds a sequence $\langle X_\alpha \mid \alpha < \omega_3\rangle$ of subsets $X_\alpha \subseteq \omega_1$ such that $X_\beta-X_\alpha$ is finite and $X_\alpha -X_\beta$ is uncountable for all $\beta < \alpha < \omega_3$. In the same way it is possible to construct along a simplified $(\omega_1,2)$-morass a ccc forcing which adds $\omega_3$ many distinct functions $f_\alpha:\omega_1 \rightarrow \omega$ such that $\{ \xi < \omega_1 \mid f_\alpha(\xi)=f_\beta(\xi)\}$ is finite for all $\alpha < \beta < \omega_3$. Another example \cite{Irrgang5} is that there exists consistently a ccc forcing of size $\omega_1$ that adds  a $0$-dimensional Hausdorff topology $\tau$ on $\omega_3$ which has spread $\omega_1$. This implies that the existence of a $0$-dimensional Hausdorff space $X$ with spread $\omega_1$ and size $2^{2^{spread(X)}}$ is consistent. 
\smallskip\\
As the examples in \cite{Irrgang4,Irrgang5} show, often consistency statements like above cannot be extended by simply raising the cardinal parameters. The reason why such a generalization could not work is that the higher-gap  case yields a higher-dimensional construction. Therefore, the finite conditions of our forcing have to fit together appropriately in more directions and that might be impossible. Hence if and how a statement generalizes to higher-gaps depends heavily on the concrete conditions. 
\smallskip\\
Morasses were introduced by R. Jensen in the early 1970's to solve the cardinal transfer problem of model theory in $L$ (see e.g. Devlin \cite{Devlin}). For the proof of the gap-2 transfer theorem a gap-1 morass is used. For higher-gap transfer theorems Jensen has developed so-called higher-gap morasses \cite{Jensen1}. In his Ph.D. thesis, the author generalized these to gaps of arbitrary size \cite{Irrgang3,Irrgang1,Irrgang2}. The theory of morasses is very far developed and very well examined. In particular it is known how to construct morasses in $L$ \cite{Devlin,Friedman,Irrgang3,Irrgang2} and how to force them \cite{Stanley2,Stanley}. Moreover, D. Velleman has defined so-called simplified morasses, along which morass constructions can be carried out very easily compared to classical morasses \cite{Velleman1984,Velleman1987a,Velleman1987b}. Their existence is equivalent to the existence of usual morasses \cite{Donder,Morgan1998}. The fact that the theory of morasses is so far developed is an advantage of the morass approach compared to historic forcing or $\rho$-functions. It allows canonical generalizations to higher cardinals, as shown below.

\section{Morasses}
As outlined in the introduction, we want to construct systems of partial orders and embeddings between them which are indexed along a ``higher-dimensional" structure. The ``two-dimensional" structure which is appropriate for our purposes is a so-called simplfied gap-1 morass. This is a very natural choice as simplified morasses can replace cardinals also in similar contexts.   Simplified morasses were introducd by D. Velleman in his papers \cite{Velleman1984,Velleman1987a} where one can also find the proofs of the following results. 
 \smallskip\\
A simplified $(\kappa, 1)$-morass is a structure $\frak{M}=\langle \langle \theta _\alpha \mid \alpha \leq \kappa \rangle , \langle \frak{F}_{\alpha\beta}\mid \alpha < \beta \leq \kappa\rangle\rangle$ satisfying the following conditions:
\smallskip\\
(P0) (a) $\theta _0=1$, $\theta _\kappa =\kappa^+$, $\forall \alpha < \kappa \ \ 0<\theta_\alpha < \kappa$.
\smallskip\\
(b) $\frak{F}_{\alpha\beta}$ is a set of order-preserving functions $f:\theta _\alpha \rightarrow \theta _\beta$.
\smallskip\\
(P1) $|\frak{F}_{\alpha\beta}| < \kappa$ for all $\alpha < \beta < \kappa$.
\smallskip\\
(P2) If $\alpha < \beta < \gamma$, then $\frak{F}_{\alpha\gamma}=\{ f \circ g \mid f \in \frak{F}_{\beta\gamma}, g \in \frak{F}_{\alpha\beta}\}$.
\smallskip\\
(P3) If $\alpha < \kappa$, then $\frak{F}_{\alpha, \alpha +1}=\{ id \upharpoonright \theta _\alpha , f_\alpha\}$ where $f_\alpha$ is such that $f_\alpha \upharpoonright \delta = id \upharpoonright \delta$ and $f_\alpha (\delta) \geq \theta _\alpha$ for some $\delta < \theta _\alpha$.
\smallskip\\
(P4) If $\alpha \leq \kappa$ is a limit ordinal, $\beta _1,\beta_2 < \alpha$ and $f_1 \in \frak{F}_{\beta _1\alpha}$, $f_2 \in \frak{F}_{\beta _2\alpha}$, then there are a $\beta _1, \beta _2 < \gamma < \alpha$, $g \in \frak{F}_{\gamma\alpha}$ and  $h_1 \in \frak{F}_{\beta _1\gamma}$, $h_2 \in \frak{F}_{\beta _2\gamma}$ such that $f_1=g \circ h_1$ and $f_2=g \circ h_2$.
\smallskip\\
(P5) For all $\alpha >0$, $\theta _\alpha =\bigcup \{ f[\theta _\beta] \mid \beta < \alpha , f \in \frak{F}_{\beta\alpha}\}$.
\bigskip\\
{\bf Lemma 2.1}
\smallskip\\
Let $\alpha < \beta \leq \kappa$, $\tau _1,\tau _2 < \theta _\alpha$, $f_1,f_2 \in \frak{F}_{\alpha \beta}$ and $f_1(\tau _1)=f_2(\tau _2)$. Then $\tau _1=\tau _2$ and $f_1\upharpoonright \tau _1 = f_2 \upharpoonright \tau _2$.
\bigskip\\
A simplified morass defines a tree $\langle T , \prec \rangle$.
\medskip\\
Let $T=\{ \langle \alpha , \nu\rangle \mid \alpha \leq \kappa , \nu < \theta _\alpha\}$.
\smallskip\\
For $t =\langle \alpha , \nu \rangle \in T$ set $\alpha (t)=\alpha$ and $\nu (t)=\nu$.
\smallskip\\
Let $\langle\alpha , \nu \rangle \prec \langle\beta , \tau\rangle$ iff $\alpha <\beta$ and $f(\nu)=\tau$ for some $f \in \frak{F}_{\alpha \beta}$.
\smallskip\\
If $s \prec t$, then $f \upharpoonright (\nu (s) +1)$ does not depend on $f$ by lemma 2.1. So we may define $\pi _{st}:=f \upharpoonright (\nu (s) +1)$.
\bigskip\\
{\bf Lemma 2.2}
\smallskip\\
The following hold:
\smallskip\\
(a) $\prec$ is a tree, $ht_T(t)=\alpha (t)$.
\smallskip\\
(b) If $t_0 \prec t_1 \prec t_2$, then $\pi _{t_0t_1}=\pi _{t_1t_2} \circ \pi _{t_0t_1}$.
\smallskip\\
(c) Let $s \prec t$ and $\pi =\pi _{st}$. If $\pi (\nu^\prime)=\tau^\prime$, $s^\prime=\langle \alpha (s), \nu ^\prime\rangle$ and $t^\prime=\langle \alpha (t), \tau^\prime\rangle$, then $s^\prime \prec t^\prime$ and $\pi _{s^\prime t^\prime}=\pi \upharpoonright (\nu^\prime +1)$.
\smallskip\\
(d) Let $\gamma \leq \kappa$, $\gamma \in Lim$. Let $t \in T_\gamma$. Then $\nu (t) +1=\bigcup \{ rng(\pi_{st})\mid s \prec t\}$.
\bigskip\\
The usual we to visualize a simplified $(\kappa,1)$-morass is as a part of a rectangle with height $\kappa+1$ and width $\kappa^+$ where every element $\langle\alpha,\nu\rangle$ of $T$ is represented by the point with coordinates $\alpha$ and $\nu$. This justifies to call a simplified gap-1 morass a two-dimensional system. Following the same idea, one can also define three-dimensional systems. This is done by D. Velleman in \cite{Velleman1987a}, where one can also find nice pictures of the two- and three-dimensional systems. 
\bigskip\\
{\bf Theorem 2.3}
\smallskip\\   
(a) If $V=L$, then there is a simplified $(\kappa ,1)$-morass for all regular $\kappa >\omega$.
\smallskip\\
(b) If $\kappa$ is an uncountable regular cardinal such that $\kappa ^+$ is not inaccessible in $L$, then there is a simplified $(\kappa ,1)$-morass.
\smallskip\\
(c) For every regular $\kappa > \omega$, there is a $\kappa$-complete (i.e. every decreasing sequence of length $<$ $\kappa$ has a lower bound) forcing $\P$ satisfying $\kappa ^+$-cc such that $\P \Vdash ($ there is a simplified $(\kappa ,1)$-morass).
\smallskip\\
{\bf Proof:} See Velleman \cite{Velleman1984,Velleman1984b}. $\Box$
\bigskip\\
For the above mentioned (three-dimensional) simplified gap-2 morasses we have analogous results. The existence of a simplified $(\kappa,2)$-morass is consistent for all regular $\kappa\geq\omega$, in particular they exist in $L$. 
\smallskip\\
Since the method of higher-dimensional forcing is related to Todorcevic's ordinal walks, the following fact is interesting.
\bigskip\\
{\bf Theorem 2.4}
\smallskip\\
If there exists a simplified $(\omega_1,1)$-morass, then $\Box_{\omega_1}$ holds.
\smallskip\\
{\bf Proof:} See Velleman \cite{Velleman}. $\Box$

\section{Two-dimensional forcing which preserves GCH}
Suppose we want to construct a ccc forcing $\P$ of size $\omega_1$. Then we can (as outlined in the introduction) proceed as follows: 
\smallskip\\
Let $\langle \sigma_{\alpha\beta}:\P_\alpha \rightarrow \P_\beta \mid \alpha < \beta < \omega_1\rangle$ be a continuous, commutative system of complete embeddings between forcings.
\smallskip\\
Let $\P$ be the direct limit of the system and assume that all $\P_\alpha$ are countable. Then every $\P_\alpha$ satifies ccc. Hence by a well-known theorem about finite support iterations $\P$ also satisfies ccc.  
\smallskip\\
Obviously, this only works for $\P$ of size $\omega_1$ because we take a direct limit of countable structures. 
\bigskip\\
What do we have to change to construct a ccc forcing $\P$ of size $\omega_2$ from countable approximations? 
\smallskip\\
1. We need some structure along which we index the approximations and which replaces the ordinal $\omega_1$ with its natural order. An appropriate structure will be given by a simplified $(\omega_1,1)$-morass. That is, we replace the linear system by a two-dimensional one.
\smallskip\\
2. We need a replacement for the continuous, commutative system of complete embeddings. We will also construct a  continuous, commutative system of  embeddings. However, not all of them will be complete.
\bigskip\\
We have to assume that there exists an $(\omega_1,1)$-morass to construct a ccc forcing $\P$ of size $\omega_2$. There exists already a very successful method to construct a ccc forcing of size $\omega_2$ assuming only $\Box_{\omega_1}$. This is Todorcevic's method of ordinal walks \cite{Todorcevic2007}. So, why is the new method interesting?
\smallskip\\
1. There exists a variant of it which guarantees that  $\P$ densely embeds into a  forcing of size $\omega_1$. Hence $\P$ preserves GCH. 
\smallskip\\
2. It has a natural generalization which allows to construct ccc forcings of size $\omega_3$. Since very little is known about possible structures on $\omega_3$, this might be interesting. 
\bigskip\\
We will first discuss the version of higher-dimensional forcing which can preserve GCH. The following is also contained in \cite{Irrgang4} where we introduced this idea.
\smallskip\\
Let $\P$ and $\Q$ be partial orders. A map $\sigma :\P \rightarrow \Q$ is called a complete embedding if
\smallskip\\
(1) $\forall  p ,p^\prime \in \P \ (p^\prime \leq p \rightarrow \sigma (p^\prime ) \leq \sigma (p))$
\smallskip\\
(2) $\forall  p ,p^\prime \in \P \ ( p$ and $p^\prime$ are incompatible $\leftrightarrow$ $\sigma (p)$ and $\sigma (p^\prime )$ are incompatible)
\smallskip\\
(3) $\forall q \in \Q \ \exists p \in \P \ \forall p^\prime \in \P \ (p^\prime \leq p \rightarrow (\sigma (p^\prime )$ and $q$ are compatible in $\Q$)).
\smallskip\\
In (3), we call $p$ a reduction of $q$ to $\P$ with respect to $\sigma$.
\medskip\\
If only (1) and (2) hold, we say that $\sigma$ is an embedding. If $\P \subseteq \Q$ such that the identity is an embedding, then we write $\P \subseteq _\bot \Q$.
\medskip\\
We say that $\P \subseteq \Q$ is completely contained in $\Q$ if $id \upharpoonright \P :\P \rightarrow \Q$ is a complete embedding.
\bigskip\\
Let $\alpha \in Lim$. A finite support (FS) iteration is a sequence $\langle \P _\xi \mid \xi \leq \alpha \rangle$ of partial orders together with a commutative system $\langle \sigma _{\xi \eta} \mid \xi < \eta \leq \alpha \rangle$ of complete embeddings $\sigma _{\xi \eta}:\P _\xi \rightarrow \P _\eta$ such that $\bigcup \{ \sigma _{\xi \eta}[\P _\xi ] \mid \xi < \eta \}=\P _\eta$ for limit $\eta$.
\medskip\\
This is the original definition by Solovay and Tennenbaum in \cite{SolovayTennenbaum}, except that they use Boolean algebras instead of partial orders. Moreover, it is well known that if $\sigma:\P_1 \rightarrow \P_2$ is a complete embedding then there is a $\P_1$-name $\dot \Q$ such that $\P_2$ and $\P_1 \ast \dot \Q$ are forcing equivalent. This leads to the more common definition of FS iterations where conditions are sequences of names. For the exact relationship between the two approaches see Kunen's textbook \cite{Kunen}, chapter VIII \S 5 and exercise K. 
\bigskip\\
An important property of FS iterations is that they preserve the $\kappa$-cc:
\medskip\\
{\bf Theorem 3.1}
\smallskip\\
Let $\langle\langle \P _\xi \mid \xi \leq \alpha \rangle ,\langle \sigma _{\xi \eta} \mid \xi < \eta \leq \alpha \rangle \rangle$ be a FS iteration. Assume that all $\P _\xi$ with $\xi < \alpha$ satisfy the $\kappa$-cc. Then $\P _\alpha$ also satisfies the $\kappa$-cc.
\smallskip\\
{\bf Proof:} See the original article by Solovay and Tennenbaum \cite{SolovayTennenbaum} or any standard textbook. $\Box$
\bigskip\\
Let $\frak{M}$ be a simplified $(\kappa ,1)$-morass. We want to define a generalization of a FS iteration which is not indexed along an ordinal but along $\frak{M}$. One way of doing this is the following definition:
\medskip\\
We say that $\langle\langle \P _\eta \mid \eta \leq \kappa ^+ \rangle ,\langle \sigma _{st} \mid s \prec t \rangle , \langle e_\alpha \mid \alpha < \kappa \rangle\rangle$ is a FS system along $\frak{M}$ if the following conditions hold:
\medskip\\
(FS1) $\langle \P _\eta \mid \eta \leq \kappa ^+  \rangle$ is a sequence of partial orders such that $\P _\eta \subseteq _\bot \P _\nu$ if $\eta \leq \nu$ and $\P _\lambda =\bigcup \{ \P _\eta\mid \eta < \lambda\}$ for $\lambda \in Lim$.
\medskip\\
(FS2) $\langle \sigma _{st} \mid s \prec t \rangle$ is a commutative system of injective embeddings $\sigma _{st}:\P _{\nu (s)+1} \rightarrow \P _{\nu (t)+1}$ such that if $t$ is a limit point in $\prec$, then $\P _{\nu (t) +1} = \bigcup \{ \sigma _{st}[\P _{\nu (s)+1}]\mid s \prec t \}$. 
\medskip\\
(FS3) $e_\alpha : \P _{\theta _{\alpha +1}} \rightarrow \P _{\theta _\alpha}$.
\medskip\\
(FS4) Let $s \prec t$ and $\pi =\pi _{st}$. If $\pi (\nu^\prime)=\tau^\prime$, $s^\prime=\langle \alpha (s), \nu ^\prime\rangle$ and $t^\prime=\langle \alpha (t), \tau^\prime\rangle$, then $\sigma _{st}:\P_{\nu (s)+1} \rightarrow \P _{\nu (t) +1}$ extends $\sigma _{s^\prime t^\prime}:\P _{\nu ^\prime +1}\rightarrow \P _{\tau^\prime +1}$.
\medskip\\
Hence for $f \in \frak{F}_{\alpha\beta}$, we may define $\sigma _f  =\bigcup \{ \sigma _{st} \mid s=\langle \alpha , \nu \rangle, t=\langle \beta , f(\nu)\rangle\}$.
\medskip\\
(FS5) If $\pi _{st}=id \upharpoonright \nu (s)+1$, then $\sigma _{st}=id \upharpoonright \P _{\nu (s) +1}$.
\medskip\\
(FS6)(a) If $\alpha < \kappa$, then $\P _{\theta_\alpha}$ is completely contained in $\P _{\theta _{\alpha +1}}$ in such a way that $e_\alpha(p)$ is a reduction of $p \in \P _{\theta _{\alpha +1}}$. 
\smallskip\\
(b) If $\alpha < \kappa $, then $\sigma _\alpha :=\sigma _{f_\alpha} :\P _{\theta _\alpha} \rightarrow \P _{\theta _{\alpha +1}}$ is a complete embedding such that $e_\alpha(p)$ is a reduction of $p \in \P _{\theta _{\alpha +1}}$.
\medskip\\
(FS7)(a) If $\alpha < \kappa$ and $p \in \P _{\theta _\alpha}$, then $e_\alpha(p)=p$.
\smallskip\\
(b) If $\alpha < \kappa$ and $p \in rng(\sigma _\alpha)$, then $e_\alpha(p)=\sigma ^{-1}_\alpha(p)$.
\bigskip\\
To simplify notation, set $\P := \P _{\kappa ^+}$.
\medskip\\
Unlike in the case of FS iterations, it is unclear how a generic extension with respect to $\P_{\kappa^+}$ can be viewed as being obtained by successive extensions. This would justify to call a FS system along $\frak{M}$ a FS {\it iteration} along $\frak{M}$.
\medskip\\
However, like in the case of FS iterations it is sometimes more convenient to represent $\P$ as a set of functions $p^\ast:\kappa \rightarrow V$ such that $p^\ast(\alpha) \in \P _{\theta _\alpha}$ for all $\alpha < \kappa$.
\bigskip\\
To define such a function $p^\ast$ from $p \in \P$ set recursively
\smallskip

$p_0=p$
\smallskip

$\nu _n(p)=min\{ \eta \mid p_n \in \P _{\eta +1}\}$
\smallskip

$t_n(p)=\langle \kappa ,\nu _n(p)\rangle$
\smallskip

$p^{(n)}(\alpha)=\sigma ^{-1}_{st}(p_n)$ if $s \in T_\alpha$, $s \prec t_n(p)$ and $p_n \in rng(\sigma_{st})$.
\medskip\\
Note that, by lemma 2.2 (a), $s$ is uniquely determined by $\alpha$ and $t_n(p)$. Hence we really define a function. Set

$\gamma _n(p)=min(dom(p^{(n)}))$.
\medskip\\
By (FS2), $\gamma _n(p)$ is a successor ordinal or $0$. Hence, if $\gamma_n(p)\neq 0$, we may define
\smallskip

$p_{n+1}=e_{\gamma_n(p)-1}(p^{(n)}(\gamma _n(p)))$.
\smallskip\\
If $\gamma_n(p)= 0$, we let $p_{n+1}$ be undefined.
\medskip\\ 
Finally, set $p^\ast=\bigcup\{ p^{(n)} \upharpoonright [\gamma _n(p) , \gamma _{n-1}(p)[ \ \mid n \in \omega \}$ where $\gamma _{-1}(p)=\kappa$.
\medskip\\
Note: If $n>0$ and $\alpha \in [\gamma _n(p),\gamma_{n-1}(p)[$, then $p^\ast(\alpha)=\sigma _{s\bar t}^{-1}(p_n)$ where $\bar t=\langle \gamma_n(p)-1,\nu _n(p)\rangle$ because $p^\ast(\alpha)=p^{(n)}(\alpha)=\sigma^{-1}_{st}(p_n)=(\sigma _{\bar tt}\circ \sigma _{s\bar t})^{-1}(p_n)=\sigma_{s\bar t}(p_n)$ where the first two equalities are just the definitions of $p^\ast$ and $p^{(n)}$. For the third equality note that $\bar t \prec t$ since $id \upharpoonright \theta _\alpha \in \frak{F}_{\alpha\beta}$ for all $\alpha < \beta \leq \kappa$ by lemma 2.3. So the equality follows from the commutativity of $\langle \sigma _{st} \mid s \prec t \rangle$. The last equality holds by (FS5).
\medskip\\
It follows from the previous observation that $\langle \gamma _n(p) \mid n \in \omega \rangle$ is decreasing. So the recursive definition above breaks down at some point, i.e. $\gamma _n(p)=0$ for some $n \in \omega$. Hence the support of $p$ which is defined by $supp(p)=\{ \gamma _n(p) \mid n \in \omega\}$ is finite.  
\bigskip\\
{\bf Theorem 3.2}
\smallskip\\
If $p^\ast(\alpha)$ and $q^\ast(\alpha)$ are compatible for $\alpha=max(supp(p) \cap supp(q))$, then $p$ and $q$ are compatible. 
\smallskip\\
{\bf Proof:} Suppose that $p$ and $q$ are incompatible. Without loss of generality let $\nu := min\{ \eta \mid p \in \P_{\eta +1}\} \leq min \{ \eta \mid q \in \P _{\eta +1} \} =:\tau$. Set $s=\langle \kappa,\nu\rangle$ and $t=\langle \kappa, \tau \rangle$. Let $t^\prime \prec t$ be minimal such that $\nu \in rng(\pi _{t^\prime t})$ and $p,q \in rng(\sigma _{t^\prime t})$. By (FS2), $t^\prime \in T_{\alpha _0 +1}$ for some $\alpha < \kappa$. Let $\pi _{t^\prime t}(\nu ^\prime)=\nu$ and $s^\prime =\langle \alpha +1, \nu ^\prime\rangle$. Let $\bar s$, $\bar t$ be the direct predecessors of $s^\prime$ and $t^\prime$ in $\prec$. Set $p^\prime = \sigma _{s^\prime s}^{-1}(p)$, $q^\prime=\sigma ^{-1}_{t^\prime t}(q)$. Then $p^\prime =p^\ast(\alpha _0+1)$, $q^\prime =q^\ast(\alpha _0+1)$ by the definition of $p^\ast$. Moreover, $p^\prime$ and $q^\prime$ are not compatible, because if $r \leq p^\prime , q^\prime$, then $\sigma _{t^\prime t}(r) \leq p,q$ by (FS2). Now, we consider several cases.
\medskip\\
{\it Case 1}: $\nu ^\prime \notin rng(\pi _{\bar t t^\prime})$
\smallskip\\
Then $\pi _{\bar s s^\prime}=id \upharpoonright \nu (\bar s)+1$ and $\sigma _{\bar s s^\prime}=id \upharpoonright \P _{\nu (\bar s)+1}$ by the minimality of $\alpha _0$. Moreover, $\bar p:=p^\prime$ and $\bar q:=e_\alpha (q^\prime)$ are not compatible, because if $r \leq p^\prime, e_\alpha (q^\prime)$, then there is $u \leq r,q^\prime,p^\prime$ by (FS6)(a). There is no difference between compatibility in $\P _{\theta _{\alpha +1}}$ and in $\P _{\nu (t^\prime)+1}$ by (FS1). Finally, note that $\bar p=p^\ast(\alpha _0)$ and $\bar q =q^\ast(\alpha _0)$ by the definition of $p^\ast$ and (FS7).
\medskip\\
{\it Case 2}:  $\nu ^\prime \in rng(\pi _{\bar t t^\prime})$ and $\pi _{\bar ss^\prime}=id \upharpoonright \nu (\bar s)+1$
\smallskip\\
Then $\pi _{\bar t t^\prime} \neq id \upharpoonright \nu (\bar t)+1$ by the minimality of $\alpha _0$ and $\bar p:=p^\prime$ and $\bar q:=e_\alpha(q^\prime)$ are not compatible (like in case 1). However, $\bar p=p^\ast(\alpha _0)$ and $\bar q =q^\ast(\alpha _0)$ by the definition of $p^\ast$ and (FS7).
\medskip\\
{\it Case 3}: $\nu ^\prime \in rng(\pi _{\bar t t^\prime})$, $\pi _{\bar ss^\prime}\neq id \upharpoonright \nu (\bar s)+1$ and $\alpha_0 +1 \notin supp(p)$
\smallskip\\
Then $\pi _{\bar t t^\prime} \neq id \upharpoonright \nu (\bar t)+1$ by the minimality of $\alpha _0$. Set $\bar p:=\sigma ^{-1}_{\bar s s^\prime}(p^\prime)$ and $\bar q = e_\alpha (q^\prime)$. Then $\bar p$ and $\bar q$ are not compatible, because if $r \leq \bar p, \bar q$, then there is $u \leq \sigma _\alpha (r),q^\prime , p^\prime$ by (FS6)(b). However, $\bar p=p^\ast(\alpha _0)$ and $\bar q =q^\ast(\alpha _0)$ by the definition of $p^\ast$ and (FS7).
\medskip\\
{\it Case 4}: $\nu ^\prime \in rng(\pi _{\bar t t^\prime})$, $\pi _{\bar ss^\prime}\neq id \upharpoonright \nu (\bar s)+1$ and $\alpha_0 +1 \notin supp(q)$ 
\smallskip\\
Then $\pi _{\bar t t^\prime} \neq id \upharpoonright \nu (\bar t)+1$. Set $\bar q:=\sigma ^{-1}_{\bar s s^\prime}(q^\prime)$ and $\bar p = e_\alpha (p^\prime)$. Then $\bar q$ and $\bar p$ are not compatible, because if $r \leq \bar p, \bar q$, then there is $u \leq \sigma _\alpha (r),p^\prime , q^\prime$ by (FS6)(b).
\medskip\\
{\it Case 5}: $\alpha_0 +1 \in supp(p) \cap supp(q)$
\smallskip\\
Then $\alpha_0 +1 = max(supp(p) \cap supp(q))$, since $\alpha_0+1 \geq max(supp(q))$ because by definition $q \in rng(\sigma _{rt})$ where $r \prec t$ and $r \in T_{max(supp(q))}$. However, $p^\prime=p^\ast(\alpha_0+1)$, $q^\prime=q^\ast(\alpha_0+1)$ are not compatible. Contradiction.
\medskip\\
So in case 5 we are finished. If we are in cases 1 - 4, we define recursively $\alpha _{n+1}$ from $p^\ast (\alpha_n)$ and $q^\ast(\alpha _n)$ in the same way as we defined $\alpha_0$ from $p$ and $q$. Like in the previous proof that $\langle \gamma _n(p)\mid n \in \omega \rangle$ is decreasing, we see that $\langle \alpha _n \mid n \in \omega \rangle$ is decreasing. Hence the recursion breaks off, we end up in case 5 and get the desired contradiction. $\Box$
\bigskip\\
{\bf Theorem 3.3}
\smallskip\\
Let $\mu ,\kappa >\omega$ be cardinals, $\kappa$ regular. Let  $\langle\langle \P _\eta \mid \eta \leq \kappa ^+ \rangle ,\langle \sigma _{st} \mid s \prec t \rangle ,\langle e_\alpha \mid \alpha < \kappa \rangle \rangle$ be a FS system along a $(\kappa ,1)$-morass $\frak{M}$. Assume that all $\P _\eta$ with $\eta < \kappa$ satisfy the $\mu$-cc. Then $\P _{\kappa ^+}$ also does.
\smallskip\\
{\bf Proof:} Let $A \subseteq \P_{\kappa^+}$ be a set of size $\mu$. Assume by the $\Delta$-system lemma that $\{ supp(p) \mid p \in A\}$ forms a $\Delta$-system with root $\Delta$. Set $\alpha = max(\Delta)$. Then $\P _{\theta _\alpha}$ satisfies the $\mu$-cc by the hypothesis of the lemma. So there are $p \neq q \in A$ such that $p^\ast(\alpha)$ and $q^\ast(\alpha)$ are compatible. Hence $p$ and $q$ are compatible by the previous lemma. $\Box$
\bigskip\\
Let $\Q =\{ p^\ast \upharpoonright supp(p) \mid p \in \P\}$.
\smallskip\\
Define a partial order $\leq$ on $\Q$ by setting $p \leq q$ iff $dom(q) \subseteq dom(p)$ and $p (\alpha) \leq q(\alpha)$ for all $\alpha \in dom(q)$.
\bigskip\\
{\bf Lemma 3.4}
\smallskip\\
Assume that $p\leq q$ implies $e_\alpha(p) \leq e_\alpha(q)$ for all $p,q \in \P_{\theta_{\alpha+1}}$ and all $\alpha \in \kappa$. Then $i:\P \rightarrow \Q, p \mapsto p^\ast \upharpoonright supp(p)$ is a dense embedding.
\smallskip\\
{\bf Proof:} By definition, $i:\P\rightarrow \Q$ is surjective. This shows density. It remains to show that it is an embedding. Assume that $p^\prime \leq p \in \P$. Then $i(p^\prime) \leq i(p)$. This is easily seen by an induction like in theorem 3.2, using the assumption that $p\leq q$ implies $e_\alpha(p) \leq e_\alpha(q)$ for all $p,q \in \P_{\theta_{\alpha+1}}$ and all $\alpha \in \kappa$. This shows (1) in the definition of embedding. Now, assume that $p$ and $p^\prime$ are compatible in $\P$. Let $q \leq p,p^\prime$. Then $i(q) \leq i(p),i(p^\prime)$, which shows one direction of the implication in (2). The other direction follows immediately from lemma 3.2. $\Box$  
\bigskip\\
{\bf Theorem 3.5}
\smallskip\\
If there exists a simplified $(\omega_1,1)$-morass, then there is a ccc forcing of size $\omega_1$ that adds an $\omega_2$-Suslin tree.
\smallskip\\
{\bf Proof:} In \cite{Irrgang4} we constructed with the method described above a ccc forcing which adds an $\omega_2$-Suslin tree. It is easily seen that the assumption of lemma 3.4 holds for this forcing. Hence the embedding $i:\P\rightarrow \Q$ is dense. $\Q$ has size $\omega_1$. Hence $\Q$ is as wanted.  $\Box$ 
\bigskip\\
It was known before that there exists a ccc forcing which adds an $\omega_2$-Suslin tree if $\Box_{\omega_1}$ holds (Todorcevic \cite{Todorcevic2007}).  
\bigskip\\
In the introduction, we also mentioned that there exists a ccc forcing which adds a Kurepa tree, if $\Box_{\omega_1}$ holds \cite{Jensen2,Jensen3,Velickovic}. If there exists a simplified $(\omega_1,1)$-morass, no forcing is needed. The morass tree $\langle T,\prec\rangle$, which we defined in section 2, is itself a Kurepa tree. 

\section{Two-dimensional forcing which destroys GCH}
Of course, many statements in whose consistency we are interested contradict GCH. Hence in proving their consistency we have to destroy GCH. So the approach we presented in section 3 is problematic because the assumption of lemma 3.4 is natural in constructions and hence difficult to avoid. In the following, we will show how to change the construction so that it is possible to destroy GCH. As example we prove
\bigskip\\
{\bf Theorem 4.1}
\smallskip\\
Assume that there exists a (simplified) $(\omega_1,1)$-morass. Then there is a ccc forcing which adds a $g:[\omega_2]^2 \rightarrow \omega$ such that $\{ \xi < \alpha \mid g(\xi,\alpha)=g(\xi,\beta)\}$ is finite for all $\alpha < \beta < \omega_2$.
\bigskip\\
This was first proved by Todorcevic using only the assumption that $\Box_{\omega_1}$ holds. He uses ordinal walks /  $\Delta$-functions.
\smallskip\\
Note that, by the Erd\"os-Rado theorem, the existence of a function $g$ like in the theorem implies the negation of $CH$.
\bigskip\\
{\bf Proof of theorem 4.1:} For $a,b \subseteq \omega_2$ let $[a,b]:=\{\langle\alpha,\gamma\rangle\mid \alpha \in a,\beta\in b ,\gamma < \alpha\}$. Set 
$$P:=\{\langle a_p,b_p,f_p \rangle \mid f_p:[a_p,b_p]\rightarrow \omega, \ a_p,b_p\subseteq\omega_2 \hbox{ finite}\}.$$
Note, that $a_p,b_p$ are not determined by $f_p$. Nevertheless, we will abuse notation and just write $p:[a_p,b_p]\rightarrow \omega$ for the condition $\langle a_p,b_p,f_p\rangle$.
\smallskip\\
We set $p \leq q$ iff $a_q \subseteq a_p$, $b_q \subseteq b_p$ and
$p(\alpha,\gamma)\neq p(\beta,\gamma) $
for all $\alpha<\beta \in a_q$ and all $\gamma \in b_p-b_q$ with $\gamma < \alpha$.
\medskip\\
It is easy to see that $P$ does not satisfy ccc. We want to thin out  $P$ to a forcing $\P$ which satisfies ccc. More precisely, we want to thin it out so that  for every $\Delta \subseteq \omega_2$
$$\P_\Delta:=\{ p \in \P \mid a_p \subseteq \Delta\}$$
satisfies ccc. Moreover, we want that there remain enough conditions that a proof like the following still works: Let $A$ be an uncountable set of conditions. Let w.l.o.g. $\{a_p \mid p \in A\}$ be a $\Delta$-system with root $\Delta$. Consider $\{ p \upharpoonright (\Delta \times \omega_2) \mid p \in A\}$. Then there are $p \neq q \in A$ such that $p \upharpoonright (\Delta \times \omega_2)$ and $q \upharpoonright (\Delta \times \omega_2)$ are compatible. Hence, $p$ and $q$ are compatible, too.
\bigskip\\
In the recursive definition of $\P$, we use the morass tree $s\prec t$ and the mappings $\pi_{st}$ to map conditions. Let more generally $\pi :\bar \theta \rightarrow \theta$ be any order-preserving map. Then $\pi :\bar\theta \rightarrow \theta$ induces maps $\pi : \bar\theta^2 \rightarrow \theta^2$ and  $\pi :\bar \theta^2  \times \omega  \rightarrow \theta^2  \times \omega$ in the obvious way:
$$\pi : \bar\theta^2 \rightarrow \theta^2 ,\quad \langle \gamma ,\delta\rangle \mapsto \langle \pi (\gamma ),\pi(\delta) \rangle$$
$$\pi :\bar \theta^2  \times \omega  \rightarrow \theta^2  \times \omega, \quad \langle x ,\epsilon\rangle \mapsto \langle \pi (x ),\epsilon \rangle.$$
\smallskip\\
We define a system  $\langle\langle \P _{\eta} \mid \eta \leq \omega_2 \rangle ,\langle \sigma _{st} \mid s \prec t \rangle \rangle$ by induction on the levels of $\langle \langle \theta_\alpha \mid \alpha \leq \omega _1\rangle,\langle \frak{F}_{\alpha\beta} \mid \alpha < \beta \leq \omega _1\rangle \rangle$ which we enumerate by $\beta \leq \omega_1$. 
\medskip\\
{\it Base Case}: $\beta =0$
\medskip\\
Then we need only to define $\P_1$.
\smallskip\\
Let $\P_1:=\{ p \in P \mid a_p, b_p \subseteq  1 \}$.  
\medskip\\
{\it Successor Case}: $\beta = \alpha +1$
\medskip\\
We first define $\P_{\theta_\beta}$. Let it be the set of all $p \in P$ such that:
\smallskip\\
(1) $a_p , b_p \subseteq \theta_\beta$
\smallskip\\
(2) $f_\alpha^{-1}[p], (id \upharpoonright\theta_\alpha)^{-1}[p] \in \P_{\theta_\alpha}$
\smallskip\\
(3) $p \upharpoonright ((\theta_\beta \setminus \theta_\alpha) \times (\theta_\alpha \setminus \delta))$ is injective
\smallskip\\
where $f_\alpha$ and $\delta$ are like in (P3) in the definition of a simplified gap-1 morass.
\bigskip\\
For $\nu \leq \theta_\alpha$, $P_\nu$ is already defined. For $\theta_\alpha <\nu \leq \theta_\beta$  set $\P_{\nu}=\{ p \in \P_{\theta_\beta} \mid a_p, b_p \subseteq \nu  \}$. 
\medskip\\
Set 
$$\sigma _{st}:\P_{\nu (s)+1} \rightarrow \P_{\nu(t)+1} , p \mapsto \pi _{st}[p].$$ \smallskip\\
{\it Limit Case}: $\beta \in Lim$
\medskip\\
For $t \in T_\beta$ set $\P_{\nu(t)+1}=\bigcup \{ \sigma _{st}[\P_{\nu (s)+1}] \mid s \prec t \}$ and $\P_\lambda =\bigcup \{ \P_\eta \mid \eta < \lambda\}$ for $\lambda \in Lim$ where $\sigma _{st}:\P _{\nu (s)+1} \rightarrow \P_{\nu(t)+1}, p \mapsto \pi _{st}[p]$.
\medskip\\
We set $\P:=\P_{\omega_2}$.
\bigskip\\
{\bf Lemma 4.2}
\smallskip\\
For $p \in P$, $p \in \P$ iff for all $\alpha <\omega_1$ and all $f\in \frak{F}_{\alpha+1,\omega_1}$
$$f^{-1}[p] \upharpoonright ((\theta_{\alpha+1} \setminus \theta_\alpha) \times (\theta_\alpha \setminus \delta_\alpha)) \hbox{ is injective}$$
where $\delta_\alpha$ is the critical point of $f_\alpha$ which is like in (P3) of the definition of a gap-1 morass.
\smallskip\\
{\bf Proof:} By induction on $\gamma \leq \omega_1$ we prove the following
\smallskip\\
{\it Claim:} $p \in \P_{\theta_\gamma}$ iff for all $\alpha < \gamma$ and all $f \in \frak{F}_{\alpha+1,\gamma}$
$$f^{-1}[p] \upharpoonright ((\theta_{\alpha+1} \setminus \theta_\alpha) \times (\theta_\alpha \setminus \delta_\alpha)) \hbox{ is injective}.$$
{\it Base case:} $\gamma=0$
\smallskip\\
Then there is nothing to prove.
\medskip\\
{\it Successor case:} $\gamma=\beta+1$
\smallskip\\
Assume first that $p \in \P_{\theta_\gamma}$. Then, by (2) in the successor step of the definition of $P_{\omega_3}$, $f^{-1}[p], (id \upharpoonright\theta_\beta)^{-1}[p] \in \P_{\theta_\beta}$. Now assume $f \in \frak{F}_{\alpha+1,\gamma}$ and $\alpha < \beta$. Then $f=f_\beta \circ f^\prime$ or $f=f^\prime$ for some $f^\prime \in \frak{F}_{\alpha+1,\beta}$ by (P2) and (P3). So by the induction hypothesis 
$$f^{-1}[p] \upharpoonright ((\theta_{\alpha+1} \setminus \theta_\alpha) \times (\theta_\alpha \setminus \delta_\alpha)) \hbox{ is injective}$$
for all $f \in \frak{F}_{\alpha+1,\gamma}$ and all $\alpha < \beta$. Moreover, if $\alpha =\beta$ then the identity is the only $f\in \frak{F}_{\alpha+1,\gamma}$. In this case 
$$f^{-1}[p] \upharpoonright ((\theta_{\alpha+1} \setminus \theta_\alpha) \times (\theta_\alpha \setminus \delta_\alpha)) \hbox{ is injective}$$
by (3) in the successor case of the definition of $\P$.
\smallskip\\
Now suppose that
$$f^{-1}[p] \upharpoonright ((\theta_{\alpha+1} \setminus \theta_\alpha) \times (\theta_\alpha \setminus \delta_\alpha)) \hbox{ is injective}$$
for all $\alpha < \gamma$ and all $f \in \frak{F}_{\alpha+1,\gamma}$. We have to prove that (2) and (3) in the successor step of the definition of $\P$ hold. (3) obviously holds by the assumption because the identity is the only function in $\frak{F}_{\gamma\gamma}=\frak{F}_{\beta+1,\gamma}$. For (2), it suffices by the induction hypothesis to show that 
$$f^{-1}[f_\beta^{-1}[p]] \upharpoonright ((\theta_{\alpha+1} \setminus \theta_\alpha) \times (\theta_\alpha \setminus \delta_\alpha)) \hbox{ is injective}$$
and
$$f^{-1}[(id\upharpoonright \theta_\beta)^{-1}[p]] \upharpoonright ((\theta_{\alpha+1} \setminus \theta_\alpha) \times (\theta_\alpha \setminus \delta_\alpha)) \hbox{ is injective}$$
for all $f \in \frak{F}_{\alpha+1,\beta}$. This, however, holds by (P2) and the assumption.
\medskip\\
{\it Limit case:} $\gamma \in Lim$
\smallskip\\
Assume first that $p \in \P_{\theta_\beta}$. Let $\alpha < \gamma$ and $f \in \frak{F}_{\alpha+1,\gamma}$. We have to prove that
$$f^{-1}[p] \upharpoonright ((\theta_{\alpha+1} \setminus \theta_\alpha) \times (\theta_\alpha \setminus \delta_\alpha)) \hbox{ is injective}.$$
By the limit step of the definition of $\P$, there are $\beta < \gamma$, $g \in \frak{F}_{\beta\gamma}$ and $\bar p\in \P_{\theta_\beta}$ such that $p=g[\bar p]$. By (P4) there are $\alpha,\beta<\delta<\gamma$, $g^\prime \in \frak{F}_{\beta\delta}$, $f^\prime\in \frak{F}_{\alpha\delta}$ and $h \in \frak{F}_{\delta\gamma}$ such that $g=h\circ g^\prime$ and $f=h\circ f^\prime$. Let $p^\prime:=g^\prime[p]$. Then, by the induction hypothesis
$$(f^\prime)^{-1}[p^\prime] \upharpoonright ((\theta_{\alpha+1} \setminus \theta_\alpha) \times (\theta_\alpha \setminus \delta_\alpha)) \hbox{ is injective}.$$
However, $(f^\prime)^{-1}[p^\prime]=(f^\prime)^{-1}[h^{-1}[p]]=f^{-1}$ and we are done. 
\smallskip\\
Now assume that
$$f^{-1}[p] \upharpoonright ((\theta_{\alpha+1} \setminus \theta_\alpha) \times (\theta_\alpha \setminus \delta_\alpha)) \hbox{ is injective}$$
for all $\alpha < \gamma$ and all $f \in \frak{F}_{\alpha+1,\gamma}$. We have to prove that $p\in \P_{\theta_\gamma}$, i.e. that there exist $\beta < \gamma$, $f\in \frak{F}_{\beta\gamma}$ and $\bar p \in \P_{\theta_\beta}$ such that $p=f[\bar p]$. However, since $p:[a_p,b_p]\rightarrow \omega$ is finite, there exist $\beta < \gamma$ and $g\in \frak{F}_{\beta\gamma}$ such that $p \in rng(f)$. Hence by the induction hypothesis it suffices to prove that $\bar p:=g^{-1}[p]\in \P_{\theta_\beta}$, i.e. that
$$f^{-1}[\bar p] \upharpoonright ((\theta_{\alpha+1} \setminus \theta_\alpha) \times (\theta_\alpha \setminus \delta_\alpha)) \hbox{ is injective}$$
for all $\alpha < \beta$ and all $f\in \frak{F}_{\alpha+1,\beta}$. So let $f \in \frak{F}_{\alpha+1,\beta}$. Then  
$$f^{-1}[\bar p] \upharpoonright ((\theta_{\alpha+1} \setminus \theta_\alpha) \times (\theta_\alpha \setminus \delta_\alpha)) =f^{-1}[g^{-1}[p]] \upharpoonright ((\theta_{\alpha+1} \setminus \theta_\alpha) \times (\theta_\alpha \setminus \delta_\alpha)) =$$
$$=(g \circ f)^{-1}[p] \upharpoonright ((\theta_{\alpha+1} \setminus \theta_\alpha) \times (\theta_\alpha \setminus \delta_\alpha))$$
which is injective by our assumption. $\Box$
\bigskip\\
For $p \in \P$ set
$$D_p=\{ \alpha < \omega_1 \mid \exists f \in \frak{F}_{\alpha+1,\omega_1} \ f^{-1}[p] \upharpoonright ((\theta_{\alpha+1} \setminus \theta_\alpha) \times (\theta_\alpha \setminus \delta_\alpha)) \neq \emptyset\}.$$
\smallskip\\
{\bf Lemma 4.3}
\smallskip\\
$D_p$ is finite for all $p \in \P$.
\smallskip\\
{\bf Proof:} For every $\langle \gamma,\xi\rangle \in dom(p)$ set $s(\gamma,\xi):=\langle \omega_1,\gamma\rangle$ and let $t(\gamma,\xi)$ be the minimal $t \prec s(\gamma,\xi)$ such that $\xi \in rng(\pi_{t,s(\gamma,\xi)})$. Then, by the inductive proof of lemma 4.2,
$$D_p:=\{ \alpha \mid \exists \langle\gamma,\xi\rangle \in dom(p) \ t(\gamma,\xi)\in T_\alpha\}.$$
Hence $D_p$ is finite because $dom(p)$ is finite. $\Box$ 
\bigskip\\
Let $\Delta \subseteq \omega_2$ be finite and $\P_\Delta=\{ p \in \P \mid a_p \subseteq \Delta\}$. We want to represent every $p \in \P_\Delta$ as a function $p^\ast: [\alpha_0,\omega_1[ \rightarrow \P$ such $p^\ast(\alpha) \in \P_{\theta_\alpha}$ for all $\alpha_0\leq \alpha < \omega_1$: Set
\smallskip

$\eta=max(\Delta)$
\smallskip

$t=\langle \omega_1,\eta\rangle$
\smallskip

$s_0=min\{ s \prec t \mid \Delta \subseteq rng(\pi_{st})\}$
\smallskip

$\alpha_0 =\alpha(s_0)$
\smallskip

$p^\ast(\alpha)=\pi^{-1}_{st}[p]$ for $\alpha _0\leq\alpha < \omega_1$
where $s \in T_\alpha$, $s \prec t$
\medskip\\
$supp(p)=$
$$\{ \alpha +1 \mid \alpha _0 \leq \alpha < \omega_1, p^\ast(\alpha +1)\neq p^\ast (\alpha), p^\ast(\alpha+1)\neq f_\alpha[p^\ast(\alpha)]\} \cup \{ \alpha_0\}$$
where $f_\alpha$ is like in (P3) of the definition of a simplified gap-1 morass.
\medskip\\
Note, that by $supp(p)$ is finite, since $p$ is finite.
\bigskip\\
{\bf Lemma 4.4}
\smallskip\\
If $p,q \in \P_\Delta$ and $p^\ast(\alpha),q^\ast(\alpha)$ are compatible in $\P_{\theta_\alpha}$ for $\alpha=max(supp(p) \cap supp(q))$, then $p$ and $q$ are compatible in $\P_\Delta$.
\smallskip\\
{\bf Proof:} Suppose $p$ and $q$ are like in the lemma, but incompatible. Let $(supp(p)\cup supp(q))-\alpha =\{ \gamma_n < \dots < \gamma_1\}$. We prove by induction on $1 \leq i \leq n$, that $p^\ast(\gamma_i)$ and $q^\ast(\gamma_i)$ are incompatible for all $1 \leq i \leq n$. Since $\gamma_n=\alpha$, this yields the desired contradiction.
\smallskip\\
Note first, that $p^\ast(\gamma_1)$ and $q^\ast(\gamma_1)$ are incompatible because otherwise $p=\pi_{st}[p^\ast(\gamma_1)]$ and $q=\pi_{st}[q^\ast(\gamma_1)]$ were incompatible (for $s \in T_{\gamma_1}$, $s \prec t$). If $\gamma_1=\alpha$, we are done. So assume that $\gamma_1 \neq \alpha$. Then either $p^\ast(\gamma_1)=\pi_{\bar ss}[p^\ast(\gamma_1-1)]$ or $q^\ast(\gamma_1)=\pi_{\bar ss}[q^\ast(\gamma_1-1)]$ where $\bar s \prec s \prec t$, $\bar s \in T_{\gamma_1-1}$ and $s \in T_{\gamma_1}$. We assume in the following that $p^\ast(\gamma_1)=\pi_{\bar ss}[p^\ast(\gamma_1-1)]$. Mutatis mutandis, the other case works the same.
\smallskip\\
{\it Claim:} \quad $p^\ast(\gamma_1-1)$ and $q^\ast(\gamma_1-1)$ are incompatible in $\P_{\theta_{\gamma_1-1}}$
\smallskip\\
Assume not. Then there is $\bar r \leq p^\ast(\gamma_1-1),q^\ast(\gamma_1-1)$ in  $\P_{\theta_{\gamma_1-1}}$ such that $a_{\bar r}=a_{p^\ast(\gamma_1-1)} \cup a_{q^\ast(\gamma_1-1)}$. Let $r^\prime:=\pi_{\bar ss}[\bar r]$. Then $r^\prime \leq \pi_{\bar ss}[p^\ast(\gamma_1-1)]=p^\ast(\gamma_1)$ and $r^\prime \leq \pi_{\bar ss}[q^\ast(\gamma_1-1)]=q^\ast(\gamma_1) \upharpoonright rng(\pi_{\bar ss})^2$. In the following we will construct an $r \leq p^\ast(\gamma_1),q^\ast(\gamma_1)$ which yields the contradiction we were looking for. Let $a_r:=a_{q^\ast(\gamma_1)} \cup a_{p^\ast(\gamma_1)}$ and $b_r:=b_{q^\ast(\gamma_1)} \cup b_{p^\ast(\gamma_1)}$. For $\langle \xi,\delta\rangle \in [a_{r^\prime} , b_{r^\prime}]$ set $r(\xi,\delta):=r^\prime(\xi,\delta)$. For $\langle \xi,\delta\rangle \in [a_{q^\ast(\gamma_1)} , b_{q^\ast(\gamma_1)}]$ set $r(\xi,\delta):=q^\ast(\gamma_1)(\xi,\delta)$. Then $r(\xi,\delta)$ is defined for all $\langle \xi,\delta\rangle \in [a_r, b_r]$ except for those in $[a_{p^\ast(\gamma_1)} -a_{q^\ast(\gamma_1)},b_{q^\ast(\gamma_1)}-rng(\pi_{\bar ss})]$. For those choose any values such that (3) in the successor step of the recursive definition of $\P$ holds. Then obviously $r\in \P_{\theta_{\gamma_1}}$. It remains to prove $r \leq p^\ast(\gamma_1),q^\ast(\gamma_1)$.  That is, we must show that
\smallskip\\ 
(1) $r(\alpha,\xi)\neq r(\beta,\xi)$ for all $\alpha<\beta \in a_{p^\ast(\gamma_1)}$ and all $\xi \in b_r-b_{p^\ast(\gamma_1)}$ with $\xi < \alpha$
\smallskip\\
(2) $r(\alpha,\xi)\neq r(\beta,\xi)$ for all $\alpha<\beta \in a_{q^\ast(\gamma_1)}$ and all $\xi \in b_r-b_{q^\ast(\gamma_1)}$ with $\xi < \alpha$.
\smallskip\\
The first statement is clear if $\xi \in b_{r^\prime}$ because $r^\prime \leq p^\ast(\gamma_1)$. So assume $\xi \notin b_{r^\prime}$. Then $\xi \notin rng(\pi_{\bar ss})$. Now, we use (P3) in the definition a simplified gap-1 morass. From (P3) and the fact that $\xi \notin rng(\pi_{\bar ss})$, $\alpha \in a_{p^\ast(\gamma_1)}$ and $\xi < \alpha$ it follows that $\pi_{\bar ss}ø\neq id \upharpoonright \nu(\bar s) +1$. Moreover, if $\delta$ is the critical point of $f_{\gamma_1-1}$ like in (P3), then $\xi \in \theta_{\gamma_1-1} \setminus \delta$ and $\alpha ,\beta \in \theta_{\gamma_1}\setminus \theta_{\gamma_1-1}$. Hence the first statement holds because of (3) in the successor step of the recursive definition of $\P$. 
\smallskip\\
The proof of the second statement is mutatis mutandis the same.  This proves the claim.
\smallskip\\
It follows from the claim, that $p^\ast(\gamma_2)$ and $q^\ast(\gamma_2)$ are incompatible. Hence we can prove the lemma by repeating this argument inductively finitely many times. $\Box$
\bigskip\\
{\bf Lemma 4.5}
\smallskip\\
$\P$ satisfies ccc.
\smallskip\\
{\bf Proof:} Let $A\subseteq \P$ be a set of size $\omega_1$. By the $\Delta$-lemma, we may assume that $\{ D_p \mid p \in A\}$ forms a $\Delta$-system with root $D$. Since for every $\alpha \in D$ there are only countably many possibilities for
$$f^{-1}[p] \upharpoonright ((\theta_{\alpha+1} \setminus \theta_\alpha) \times (\theta_\alpha \setminus \delta_\alpha)),$$
we may moreover assume that for all $\alpha\in D$, all $f \in \frak{F}_{\alpha+1,\omega_1}$ and all $p,q \in A$
$$f^{-1}[p] \upharpoonright ((\theta_{\alpha+1} \setminus \theta_\alpha) \times (\theta_\alpha \setminus \delta_\alpha)) =f^{-1}[q] \upharpoonright ((\theta_{\alpha+1} \setminus \theta_\alpha) \times (\theta_\alpha \setminus \delta_\alpha)).$$
By the $\Delta$-system lemma, we may assume that $\{ a_p \mid p \in A\} \subseteq \omega_2$ forms a $\Delta$-system with root $\Delta_1$. Consider $A^\prime:=\{p \upharpoonright (\Delta_1 \times \omega_2)\mid p \in A\}$. By the $\Delta$-system lemma we may also assume that $\{ supp(p) \mid p \in A^\prime\} \subseteq \omega_1$ forms a $\Delta$-system with root $\Delta_2$. Let $\alpha=max(\Delta_2)$. Since $\P_{\theta_\alpha}$ is countable, there are $q_1 \neq q_2\in A^\prime$ such that $q^\ast_1(\alpha)=q_2^\ast(\alpha)$. Hence $q_1\neq q_2 \in A^\prime$ are compatible by lemma 5.3. Assume that $q_1=p^\ast_1 \upharpoonright \Delta_2$ and $q_2 = p^\ast_2 \upharpoonright \Delta_2$ with $p_1,p_2 \in A$. We can define $p \leq p_1,p_2$ as follows: $a_p=a_{p_1}\cup a_{p_2}$, $b_p=b_{p_1}\cup b_{p_2}$, $p \upharpoonright (a_{p_1}\times b_{p_1})=p_1$,   $p \upharpoonright (a_{p_2}\times b_{p_2})=p_2$. We still need to define $p$ on $[a_p, b_p] -((a_{p_1} \times b_{p_1}) \cup (a_{p_2}\times b_{p_2}))$. We do this in such a way that the new values are different from the old ones and distinct among each other.
\smallskip\\
Obviously $p \leq p_1,p_2$. It remains to prove that $p\in \P$. For this we use lemma 4.2. That is, we have to show that for all $\alpha <\omega_1$ and all $f\in \frak{F}_{\alpha+1,\omega_1}$
$$f^{-1}[p] \upharpoonright ((\theta_{\alpha+1} \setminus \theta_\alpha) \times (\theta_\alpha \setminus \delta_\alpha)) \hbox{ is injective}.$$
However, if $\alpha \in D$, then this holds by our first thinning-out of $A$. If $\alpha \notin D$, then it holds because of the way in which we defined $p$ on $[a_p, b_p] -((a_{p_1} \times b_{p_1}) \cup (a_{p_2}\times b_{p_2}))$.
$\Box$
\bigskip\\
{\bf Lemma 4.6}
\smallskip\\
Let $p \in \P$ and $\alpha,\beta \in \omega_2$. Then there exists $q \leq p$ such that $\alpha \in a_q$ and $\beta \in b_q$.
\smallskip\\
{\bf Proof:}  Let $a_q=a_p \cup \{ \alpha\}$, $b_q=b_p\cup \{ \beta\}$ and $q \upharpoonright (a_p \times b_p) =p$. We have to define $q(\alpha, \beta)$ on $[a_q, b_q] -((a_{p_1} \times b_{p_1}) \cup (a_{p_2}\times b_{p_2}))$. We do this in such a way that the new values are different from the old ones and distinct among each other. Then obviously $q \leq p$ and $q \in \P$ by lemma 4.2.  $\Box$  
\bigskip\\
Now, we are ready to prove theorem 4.1. Of course, $\P$ is the forcing which we defined above. Let $G$ be $\P$-generic and set $f=\bigcup\{ p \mid p \in G\}$. By lemma 4.5, cardinals are preserved. By lemma 4.6, $f$ is defined on all of $[\omega_2]^2$. By the definition of $\leq$, $f$ is as wanted. $\Box$ (Theorem 4.1)
\bigskip\\
{\bf Theorem 4.7}
\smallskip\\
There exists consistently a chain $\langle X_\alpha \mid \alpha < \omega_2\rangle$ such that $X_\alpha \subseteq \omega_1$, $X_\beta -X_\alpha$ is finite and $X_\alpha - X_\beta$ has size $\omega_1$ for all $\beta < \alpha < \omega_2$.
\smallskip\\
{\bf Proof:} For a proof with higher-dimensional forcing see \cite{Irrgang5} and \cite{Irrgang6}. $\Box$
\bigskip\\
The theorem was first proved by Koszmider \cite{Koszmider} using ordinal walks. Note, that if we set $Y_\alpha =X_{\alpha+1}-X_\alpha$, then $\{ Y_\alpha \mid \alpha \in \omega_2\}$ forms a family of $\omega_2$-many uncountable subsets of $\omega_1$ such that $Y_\alpha \cap Y_\beta$ is finite whenever $\alpha \neq \beta \in \omega_2$. Hence, by a result of Baumgartner's \cite{Baumgartner}, GCH cannot hold. 
\bigskip\\
{\bf Theorem 4.8}
\smallskip\\
There exists consistently a family $\langle f_\alpha \mid \alpha < \omega_2\rangle$ of functions such that $\{ \xi < \omega_1 \mid f_\alpha(\xi)=f_\beta(\xi)\}$ is finite for all $\alpha < \beta < \omega_2$.
\smallskip\\
{\bf Proof:} For a proof with higher-dimensional forcing see \cite{Irrgang6}. $\Box$
\bigskip\\ 
It is known that there can be families $\{ f_\alpha:\omega_1 \rightarrow \omega \mid \alpha \in \kappa\}$ of arbitrary prescribed size such that $\{ \xi < \omega_1 \mid f_\alpha(\xi)=f_\beta(\xi)\}$ is finite for all $\alpha < \beta < \kappa$. This was proved by J. Zapletal \cite{Zapletal} using proper forcing and Todorcevic's method of side conditions \cite{Todorcevic1985,Todorcevic1989}. Moreover, a result of Baumgartner's \cite{Baumgartner} is that for any given $\kappa$ there is consistently a family of size $\kappa$ of cofinal subsets of $\omega_1$ with pairwise finite intersections. Obviously, our family is a family of $\omega_2$-many uncountable subsets of $\omega_1\times \omega$ such that the intersection of two distinct members is always finite. Hence, as mentioned already above, by a result of Baumgartner's \cite{Baumgartner} GCH cannot hold.  
\smallskip\\
However, not everything can be done by ccc forcings. For example, Koszmider proved that if CH holds, then there is no ccc forcing that adds a sequence of $\omega_2$ many functions $f:\omega_1 \rightarrow \omega_1$ which is ordered by strict domination mod finite. However, he is able to produce a proper forcing which adds such a sequence \cite{Koszmider2000} using his method of side conditions in morasses which is an extension of Todorcevic's method of side conditions in models. More on the method can be found in Morgan's paper \cite{Morgan2006}. In the context of our approach, this raises the question if it is possible to define something like a countable support iteration along a morass. Finally, we should mention that by a theorem of Shelah's \cite{Shelah}, there cannot be  a sequence of $\omega_3$ many functions $f:\omega_1 \rightarrow \omega_1$ which is ordered by strict domination mod finite. 
\bigskip\\
{\bf Theorem 4.9}
\smallskip\\
There exists consistently an $(\omega,\omega_2)$-superatomic Boolean algebra.
\bigskip\\
This was first proved by Baumgartner and Shelah \cite{BaumgartnerShelah} independently of ordinal walks. Later Todorcevic \cite{Todorcevic2007} found in the presence of $\Box_{\omega_1}$ a ccc forcing which proves the consistency.
\smallskip\\
A superatomic Boolean algebra (abbreviated sBa) is a Boolean algebra in which every subalgebra is atomic. A Boolean algebra $B$ is superatomic iff its $\alpha$-th Cantor-Bendixon ideal $I_\alpha = B$ for some $\alpha$. The Cantor-Bendixon ideals $I_\alpha$ are defined by induction on $\alpha$ as follows: $I_0=\{ 0 \}$; if $\alpha = \beta +1$, then $I_\alpha$ is the ideal generated by $I_\beta \cup \{ b \in B \mid b/I_\beta$ is an atom in $B/I_\beta \}$; and if $\alpha \in Lim$, then $I_\alpha = \bigcup \{ I_\beta \mid \beta < \alpha \}$. This corresponds to the Cantor-Bendixon derivative of topological spaces. In fact, a Boolean algebra $B$ is superatomic iff its Stone space $S(B)$ is scattered. So results on superatomic Boolean algebras transfer readily to scattered spaces.
\smallskip\\
The height $ht(B)$ of a sBa is the least $\alpha$ such that the set of atoms of $B/I_\alpha$ is finite (or equivalently $B=I_{\alpha +1}$). For every $\alpha < ht(B)$ let $wd_\alpha (B)$ be the cardinality of the set of atoms in $B/I_\alpha$. If $\alpha$ is an ordinal with $\alpha \neq 0$, we say that $B$ is a $(\kappa , \alpha )$-sBa, if $ht(B)=\alpha$ and $wd_\beta (B) \leq \kappa$ for all $\beta < \alpha$. An $(\omega ,\omega _1)$-sBa is called thin tall. An $(\omega ,\omega _2)$-sBa is called thin-very tall.
\smallskip\\
The question of existence of a thin-tall sBas was posed by Telgarsky (1968). It was independently answered positively by Rajagopalan \cite{Rajagopalan} and by Juhasz and Weiss \cite{JuhaszWeiss}. They even proved that there exists an $(\omega , \alpha )$-sBa for each $\alpha < \omega _2$. Since $wd_\alpha (B) =\omega$ implies $card(B/I_\alpha )\leq 2^\omega$ \cite{Roitman}, this is the best we can hope for in $ZFC$. In \cite{KoepkeMartinez}, Koepke and Martinez showed that the existence of a $(\kappa ,\kappa ^+)$-sBa follows from the existence of a simplified $(\kappa ,1)$-morass. Martinez \cite{Martinez} showed that it is consistent that for all $\alpha < \omega_3$ there is an $(\omega,\alpha)$-sBa. Like long increasing chains modulo flat ideals \cite{Shelah}, also superatomic Boolean algebras are related to pcf theory \cite{JechShelah,Foreman}.
\bigskip\\
We now come to the {\bf proof of theorem 4.9}. Assume that there is a partial order $<_B$ on $\omega_2$ such that: 
\smallskip\\
(a) $\forall \alpha,\beta \in \omega_2 \ \alpha <_B\beta \rightarrow \alpha < \beta$.
\smallskip\\
(b) If $\alpha,\beta \in [\omega\gamma,\omega\gamma+\omega)$, then $\alpha \not <_B \beta$.
\smallskip\\
(c) Whenever $\alpha,\beta$ are compatible in $<_B$, then they have an infimum $i(\alpha,\beta)$.
\smallskip\\
(d) If $\alpha \in [\omega\gamma,\omega\gamma+\omega)$, then there exist for all $\delta < \gamma$ infinitely many  $\beta <_B \alpha$ such that $\beta \in [\omega\delta,\omega\delta+\omega)$. 
\bigskip\\
It is well-known \cite{KoepkeMartinez,Foreman,BaumgartnerShelah}, that if there exists such a partial order, then there also is an $(\omega,\omega_2)$-sBa.
\bigskip\\
Let $P$ be the set of all finite (strict) partial orders $p=\langle x_p,<_p\rangle$ such that:
\smallskip\\
(a) $\forall \alpha,\beta \in \omega_2 \ \alpha <_p\beta \rightarrow \alpha < \beta$.
\smallskip\\
(b) If $\alpha,\beta \in [\omega\gamma,\omega\gamma+\omega)$, then $\alpha \not <_p \beta$.
\smallskip\\
(c) Whenever $\alpha,\beta$ are compatible in $<_p$, then they have an infimum $i_p(\alpha,\beta)$.
\medskip\\
For $p,q\in P$, we set $p\leq q$ iff
\smallskip\\
(i) $x_q \subseteq x_p$, $<_p \upharpoonright x_q=<_q$
\smallskip\\
(ii) If $\alpha,\beta\in x_q$ are compatible in $<_p$, then they are compatible in $<_q$ and $i_p(\alpha,\beta)=i_q(\alpha,\beta)$.
\bigskip\\
$P$ is obviously the natural forcing to add a partial order $<_B$ as described above. Like in all our other examples, this forcing does not satisfy ccc. Therefore, we proceed like in the proof of theorem 4.1. The similarity between the two forcings is even more obvious, if we consider the functions $f_p:[a_p,b_p] \rightarrow 2$ instead of $p$ as conditions, where $f_p$ is defined as follows. We set $a_p=\{ \beta \mid \exists \alpha \ \alpha <_p \beta\}$, $b_p=\{ \alpha \mid \exists \beta \ \alpha <_p \beta\}$, and $f_p(\alpha,\beta)=1$ iff $\alpha <_p\beta$.   
\smallskip\\
As in the proof of theorem 4.1, we can map the $f_p$'s along $\pi$.  We define a system  $\langle\langle P _{\eta} \mid \eta \leq \omega_3 \rangle ,\langle \sigma _{st} \mid s \prec t \rangle \rangle$ by induction on the levels of $\langle \langle \theta_\zeta \mid \zeta \leq \omega _1\rangle,\langle \frak{G}_{\zeta\xi} \mid \zeta < \xi \leq \omega _2\rangle \rangle$ which we enumerate by $\beta \leq \omega_2$. 
\medskip\\
{\it Base Case}: $\beta =0$
\medskip\\
Then we need only to define $P_1$.
\smallskip\\
Let $P_1:=\{ p \in P \mid a_p, b_p \subseteq  1 \}$.  
\medskip\\
{\it Successor Case}: $\beta = \alpha +1$
\medskip\\
We first define $P_{\theta_\beta}$. Let it be the set of all $p \in P$ such that:
\smallskip\\
(1) $a_p , b_p \subseteq \theta_\beta$
\smallskip\\
(2) $f_\alpha^{-1}[p], (id \upharpoonright\theta_\alpha)^{-1}[p] \in P_{\theta_\alpha}$
\smallskip\\
(3) $\forall \alpha \in \theta_\alpha \setminus \delta \ card(\{ \beta \in \theta_\beta \setminus \theta_\alpha \mid f_p(\alpha,\beta)=1\}) \leq 1$
\smallskip\\
where $f_\alpha$ and $\delta$ are like in (P3) in the definition of a simplified gap-1 morass.
\bigskip\\
For $\nu \leq \theta_\alpha$, $P_\nu$ is already defined. For $\theta_\alpha <\nu \leq \theta_\beta$  set $P_{\nu}=\{ p \in P_{\theta_\beta} \mid a_p, b_p \subseteq \nu  \}$. 
\medskip\\
Set 
$$\sigma _{st}:P_{\nu (s)+1} \rightarrow P_{\nu(t)+1} , p \mapsto \pi _{st}[p].$$ \smallskip\\
{\it Limit Case}: $\beta \in Lim$
\medskip\\
For $t \in T_\beta$ set $P_{\nu(t)+1}=\bigcup \{ \sigma _{st}[P_{\nu (s)+1}] \mid s \prec t \}$ and $P_\lambda =\bigcup \{ P_\eta \mid \eta < \lambda\}$ for $\lambda \in Lim$ where $\sigma _{st}:P _{\nu (s)+1} \rightarrow P_{\nu(t)+1}, p \mapsto \pi _{st}[p]$.
\bigskip\\
The rest of the proof is now like the proof of theorem 4.1. $\Box$

\section{Three-dimensional forcing}

Let $X$ be a topological space. Its spread is defined by 
$$spread(X)=sup\{ card(D) \mid D \hbox{ discrete subspace of }X\}.$$
\smallskip\\
{\bf Theorem 5.1} (Hajnal,Juhasz - 1967)
\smallskip\\
If $X$ is a Hausdorff space, then $card(X) \leq 2^{2^{spread(X)}}$.
\bigskip\\  
In his book "Cardinal functions in topology" \cite{Juhasz}, Juhasz explicitly asks if the second exponentiation is really necessary. This was answered by Fedorcuk \cite{Fedorcuk}.
\bigskip\\
{\bf Theorem 5.2}
\smallskip\\
In $L$, there exists a $0$-dimensional Hausdorff (and hence regular) space with spread $\omega$ of size $\omega_2=2^{2^{spread(X)}}$.
\bigskip\\
This is a consequence of $\diamondsuit$ (and GCH).  
\bigskip\\
There was no such example for the case $spread(X)=\omega_1$. Three-dimensional forcing yields the following:
\bigskip\\
{\bf Theorem 5.3}
\smallskip\\
If there is a simplified $(\omega_1,2)$-morass, then there exists a ccc forcing of size $\omega_1$ which adds a $0$-dimensional Hausdorff space $X$ of size $\omega_3$ with spread $\omega_1$.
\smallskip\\
{\bf Proof:} See \cite{Irrgang5}. $\Box$
\bigskip\\
Hence there exists such a forcing in $L$. By the usual argument for Cohen forcing, it preserves $GCH$. So the existence of a $0$-dimensional Hausdorff space with spread $\omega_1$ and size $2^{2^{spread(X)}}$ is consistent.
\bigskip\\
{\bf Theorem 5.4}
\smallskip\\
It is consistent that there exists a sequence $\langle X_\alpha \mid \alpha < \omega_3\rangle$ of subsets $X_\alpha \subseteq \omega_1$ such that $X_\beta-X_\alpha$ is finite and $X_\alpha -X_\beta$ is uncountable for all $\beta < \alpha < \omega_3$.
\smallskip\\
{\bf Proof:} See \cite{Irrgang6}. $\Box$
\bigskip\\
P. Koszmider proved with the method of ordinal walks that there exists consistently a chain $\langle X_\alpha \mid \alpha < \omega_2\rangle$ such that $X_\alpha \subseteq \omega_1$, $X_\beta -X_\alpha$ is finite and $X_\alpha - X_\beta$ has size $\omega_1$ for all $\beta < \alpha < \omega_2$.
\bigskip\\
{\bf Theorem 5.5}
\smallskip\\
Assume that there exists a simplified $(\omega_1,2)$-morass. Then there is a ccc forcing which adds $\omega_3$ many distinct functions $f_\alpha:\omega_1 \rightarrow \omega$ such that $\{ \xi < \omega_1 \mid f_\alpha(\xi)=f_\beta(\xi)\}$ is finite for all $\alpha < \beta < \omega_3$.
\smallskip\\
{\bf Proof:} See \cite{Irrgang6}. $\Box$

\section{Open problems}
We write
$\kappa \rightarrow (\sigma:\tau)^2_\gamma$
for: Every partition $f:[\kappa]^2\rightarrow \gamma$ has a homogeneous set $[A;B]:=\{ \{\alpha,\beta\}\mid \alpha \in A,\beta\in B\}$ where $\alpha < \beta$ for all $\alpha \in A$ and $\beta \in B$ and $card(A)=\sigma$, $card(B)=\tau$, i.e. $f$ is constant on $[A;B]$.
\smallskip\\
We write $\kappa \not\rightarrow (\sigma:\tau)^2_\gamma$ for the negation of this statement.
\smallskip\\
A classical result concerned with our partition relation is the consistency of $\omega_2\not\rightarrow(\omega_1:2)^2_\omega$ and $\omega_3\not\rightarrow(\omega_2:2)^2_{\omega_1}$ (even with GCH). This was proved by A. Hajnal by forcing (cf. \cite{HajnalLarson}, theorem 4.3 and \cite{Unsolved2}). In view of the methods we apply, it might be interesting to know that $\omega_2\not\rightarrow (\omega_1:\omega)^2_2$ and $\omega_3\not\rightarrow (\omega_2:\omega_1)^2_2$ hold in $L$. This was shown by Rebholz \cite{Rebholz} using morasses and diamond. \bigskip\\
{\bf Question 6.1}
\smallskip\\
Assume that there is an $(\omega_1,2)$-morass. Does then exist a ccc forcing with finite conditions which forces $\omega_3\not\rightarrow (\omega :2)^2_{\omega_1}$? 
\bigskip\\
Note, that by theorem 4.1 there is a ccc forcing for $\omega_2\not \rightarrow (\omega:2)^2_\omega$ if there exists an $(\omega_1,1)$-morass. In \cite{Zapletal}, Zapletal states that the existence of a family $\{ f_\alpha:\omega_2 \rightarrow \omega_1 \mid \alpha \in \omega_3\}$ such that $\{ \xi < \omega_1 \mid f_\alpha(\xi)=f_\beta(\xi)\}$ is finite for all $\alpha \neq \beta < \omega_3$ is consistent. Using $\Box_{\omega_2}$, this implies the consistency of $\omega_3 \not\rightarrow (\omega : 2)^2_{\omega_1}$. The author thanks H.-D. Donder for pointing this out to him. However, Zapletal's forcing is not ccc. Finally, we should mention that $\omega_3 \not \rightarrow (\omega:2)^2_\omega$ is easily seen to be inconsistent.  
\bigskip\\
Another generalisation of $\omega_2\not \rightarrow (\omega:2)^2_\omega$ leads to the following question \cite{Todorcevic1991,Todorcevic1998,Todorcevic2007}. 
\bigskip\\
{\bf Question 6.2} (Todorcevic)
\smallskip\\
Can there exist (consistently) a function $f:\omega_3 \times \omega_3\rightarrow \omega$ such that $f$ is not constant on any rectangle $A\times B$ with infinite $A,B \subseteq \omega_3$?
\bigskip\\
After the discussion of the possible heights and widths of superatomic Boolean algebras in section 4, the following question remains open. This is a quite famous open problem which has its own chapter by J. Bagaria \cite{ProblemsTopology}  in the book ``Open problems in topology II''.
\bigskip\\  
{\bf Question 6.3} (famous)
\smallskip\\
Is the existence of an $(\omega,\omega_3)$-superatomic or an $(\omega_1,\omega_3)$-superatomic Boolean algebra consistent?
\bigskip\\
Even though the method generalizes straightforwardly to higher-dimensions, this is not true for the consistency statements. The reason is that the conditions of the forcing have to fit together in more directions, if we go to higher dimensions.  
\smallskip\\
{\bf Question 6.4} Can we find good applications in dimensions higher than two (three)?
\bigskip\\
Even though the method is inspired by iterated forcing, all my examples use basically finite sets of ordinals as conditions. Hence we do not need names for conditions.
\smallskip\\
{\bf Question 6.5} Can we find an application which uses names for conditions?
 \bigskip\\
 As we pointed out in section 3, there are consistency results that cannot be proved with ccc forcing.
 \bigskip\\
 {\bf Question 6.6} Can we do higher-dimensional forcing with countable support?  

\bibliographystyle{plain}
\bibliography{biblio}
\end{document}